\documentclass[12pt]{article}%
\usepackage{amsmath,amssymb}
\usepackage[applemac]{inputenc}
\usepackage{amsmath,amssymb,fullpage}
\usepackage{color}
\usepackage{amsmath}
\usepackage{amsfonts}
\usepackage{amssymb}
\usepackage{graphicx}%
\providecommand{\U}[1]{\protect\rule{.1in}{.1in}}

\newtheorem{definition}{Definition}[section]

\newtheorem{theorem}[definition]{Theorem}
\newtheorem{problem}[definition]{Problem}
\newtheorem{remark}[definition]{ \it Remark}

\newtheorem{proposition}[definition]{Proposition}

\numberwithin{equation}{section}

\def\1B{\text{1\!\!I}}

\begin{document}

\title{New approach to optimal control of stochastic Volterra integral equations }
\author{Nacira Agram$^{1,2,3}$, Bernt Øksendal$^{1,3}$ and Samia Yakhlef$^{2}$ }
\date{6 December 2018
\vskip 0.2cm
Final version will be published in Stochastics}
\maketitle

\begin{abstract}
We study optimal control of stochastic Volterra integral equations (SVIE) with
jumps by using Hida-Malliavin calculus. \newline- We give conditions under
which there exists unique solutions of such equations.\newline- Then we prove
both a sufficient maximum principle (a verification theorem) and a necessary
maximum principle via Hida-Malliavin calculus.\newline- As an application we
solve a problem of optimal consumption from a cash flow modelled by an SVIE.

\end{abstract}

\footnotetext[1]{Department of Mathematics, University of Oslo, P.O. Box 1053
Blindern, N--0316 Oslo, Norway. Email: \texttt{naciraa@math.uio.no,
oksendal@math.uio.no.}}

\footnotetext[2]{University Mohamed Khider of Biskra, Algeria. Email:
\texttt{samiayakhelef@yahoo.fr.}}

\footnotetext[3]{This research was carried out with support of the Norwegian
Research Council, within the research project Challenges in Stochastic
Control, Information and Applications (STOCONINF), project number 250768/F20.}



\section{Introduction}

Stochastic Volterra integral equations (SVIE) are a special type of integral
equations. They represent interesting models for stochastic dynamics with
memory, with applications to e.g. engineering, biology and finance.\newline In
this work, we consider the problem of optimal control of stochastic Volterra
integral equations of the form%
\[%
\begin{array}
[c]{cc}%
X^{u}(t) & =\xi(t)+%
{\textstyle\int_{0}^{t}}
b\left(  t,s,X^{u}(s),u(s)\right)  ds+%
{\textstyle\int_{0}^{t}}
\sigma\left(  t,s,X^{u}(s),u(s)\right)  dB(s)\\
& +%
{\textstyle\int_{0}^{t}}
{\textstyle\int_{\mathbb{R}_{0}}}
\gamma\left(  t,s,X^{u}(s^-),u(s^-),\zeta\right)  \tilde{N}(ds,d\zeta);t\in
\lbrack0,T],
\end{array}
\]
where $T>0$ is a fixed given terminal time, and the process $u(t)$ is our
control process. Here $B(t)$ and $\tilde{N}(dt,d\zeta):=N(dt,d\zeta
)-\nu(d\zeta)dt$ is a Brownian motion and an independent compensated Poisson
random measure, respectively, jointly defined on a filtered probability space
$(\Omega,\mathcal{F},\mathbb{F}=\{\mathcal{F}_{t}\}_{t\geq0},P)$ satisfying
the usual conditions. The measure $\nu$ is the Lévy measure of the jump
measure $N$.\\
For notational convenience, in the following we will always assume that in the integrals with
respect to $\tilde{N}(ds,d\zeta)$ the predictable version of the integrand is used.
\newline \newline
The problem is to find a control $\hat{u}$ which
maximises the performance functional $J(u)$ defined by
\[%
\begin{array}
[c]{ll}%
J(u) & =\mathbb{E[}{%
{\textstyle\int_{0}^{T}}
}\text{ }f(t,X^{u}(t),u(t))dt+h(X^{u}(T))].
\end{array}
\]
By using the maximum principle, we obtain an adjoint equation which is a
backward equation of Volterra type. In general a backward stochastic Volterra
integral equation (BSVIE) in the unknown trippel $(p(t),q(t,s),r(t,s,\cdot
));0\leq t\leq s\leq T$ has the form%
\[%
\begin{array}
[c]{cc}%
p(t) & =F(t)+%
{\textstyle\int_{t}^{T}}
g(t,s,p(s),q(t,s),r(t,s,\cdot))ds-%
{\textstyle\int_{t}^{T}}
q(t,s)dB(s)\\
& -%
{\textstyle\int_{t}^{T}}
{\textstyle\int_{\mathbb{R}_{0}}}
r(t,s^-,\zeta)\tilde{N}(ds,d\zeta);t\in\lbrack0,T].
\end{array}
\]
 
An equivalent formulation of this is that%
\[%
\begin{array}
[c]{cc}%
p(t) & =\mathbb{E[}F(t)+%
{\textstyle\int_{t}^{T}}
g(t,s,p(s),q(t,s),r(t,s,\cdot))ds|\mathcal{F}_{t}].
\end{array}
\]
One of the many interesting motivations of such an equation is the
\emph{recursive utility}: \newline For a consumption process $c(t)\geq0$, we
consider its Volterra recursive utility process $p(t)$ defined by
\[%
\begin{array}
[c]{cc}%
p(t) & =\mathbb{E[}F(t)+%
{\textstyle\int_{t}^{T}}
\varphi(t,s,p(s),c(s))ds|\mathcal{F}_{t}],
\end{array}
\]
where $F$ and $\varphi$ are given functions. This is an extension of the
classical recursive utility concept of Duffie and Epstein \cite{DE} to
Volterra integral equation with jumps. Finding an optimal consumption rate
$\hat{c}$ which maximises the total Volterra recursive utility
\[
U(\hat{c})=p(0),
\]
is an interesting problem in mathematical finance.\newline The problem of
optimal control of SVIE has been studied by several authors. See e.g. Yong
\cite{Yong1},\cite{Yong 2}, Agram \textit{et al} \cite{AO},\cite{AOY}. In
contrast to Yong \cite{Yong1} and \cite{Yong 2}, we do obtain a sufficient and
a necessary maximum principle of the classical type.\newline The paper which
is closest in content to the current paper is Agram and Øksendal \cite{AO}.
But the current paper differs from \cite{AO} in an essential way: In the
current paper both the Hamiltonian and the associated adjoint equation are
simpler and much easier to deal with, because neither of them involve any
Hida-Malliavin derivatives. Thus, the Hida-Malliavin derivatives are only used
in the proofs and do not appear in the final formulation of the maximum
principles. Moreover,
the adjoint equation obtained
\color{black}
 here is a standard BSVIE. The price we have to pay for this simplification, is
that we have to assume smoothness with respect to $t$ of the two components
$q(t,s),r(t,s,\zeta)$ of the solution $(p,q,r)$ of the associated (linear)
BSVIE. It is not clear to what extent smoothness properties hold for the
solutions of BSVIEs in general. However, in the paper by Hu and Øksendal
\cite{ho}, it is shown that the required smoothness holds for linear BSVIEs,
under certain conditions.\newline We outline the content of this
paper:\newline In Section 2 we give some background about Hida-Malliavin
calculus. \newline In Section 3 we first give conditions under which there
exists a unique solution of such an SVIE. Then we prove both a sufficient
maximum principle (a verification theorem) and a necessary maximum principle
via Hida-Malliavin calculus. \newline Finally, in Section 4 we illustrate our
results by solving a problem of optimal consumption from a cash flow modelled
by an SVIE.

\section{Framework}

Throughout this work, we will use the following spaces:

\begin{itemize}
\item $\mathcal{S}^{2}$ is the set of ${\mathbb{R}}$-valued $\mathbb{F}%
$-adapted càdlàg processes $(X(t))_{t\in\lbrack0,T]}$ such that
\[
{\Vert X\Vert}_{\mathcal{S}^{2}}^{2}:={\mathbb{E}}[\sup_{t\in\lbrack
0,T]}|X(t)|^{2}]~<~\infty\;.
\]

\item $\mathbb{L}^{2}$ is the set of ${\mathbb{R}}$-valued $\mathbb{F}%
$-adapted processes $\{ Q(t,s)\}_{(t,s)\in\lbrack0,T]^{2}}$ such that
\[
\Vert Q\Vert_{\mathbb{L}^{2}}^{2}:={\mathbb{E}}[%
{\textstyle\int_{0}^{T}}
{\textstyle\int_{t}^{T}}
|Q(t,s)|^{2}dsdt]<~\infty\;.
\]
{}

\item $\mathbb{L}_{\nu}^{2}$ is the set of Borel functions $K:%
\mathbb{R}
_{0}\rightarrow%
\mathbb{R}
,$ such that
\[
\parallel K\parallel_{\mathbb{L}_{\nu}^{2}}^{2}:=%
{\textstyle\int_{\mathbb{R}_{0}}}
K(\zeta)^{2}\nu(d\zeta)<\infty,
\]
where $%
\mathbb{R}
_{0}:=%
\mathbb{R}
\setminus\{0\}.$

\item $\mathbb{H}_{\nu}^{2}$ is the set of $\mathbb{F}$-adapted predictable
processes$\ R:[0,T]^{2}\times%
\mathbb{R}
_{0}\times\Omega\rightarrow\mathbb{R},$ such that $\mathbb{E[}%
{\textstyle\int_{0}^{T}}
{\textstyle\int_{t}^{T}}
{\textstyle\int_{\mathbb{R}_{0}}}
|R(t,s,\zeta)|^{2}\nu(d\zeta)dsdt]<\infty.$ We equip $\mathbb{H}_{\nu}^{2}$
with the norm
\[
\parallel R\parallel_{\mathbb{H}_{\nu}^{2}}^{2}:=\mathbb{E[}%
{\textstyle\int_{0}^{T}}
{\textstyle\int_{t}^{T}}
{\textstyle\int_{\mathbb{R}_{0}}}
|R(t,s,\zeta)|^{2}\nu(d\zeta)dsdt].
\]

\end{itemize}

\subsection{The generalized Hida-Malliavin derivative}

The Malliavin derivative $D_{t}$ was originally introduced by Malliavin
\cite{M} as a stochastic calculus of variation used to prove results about
smoothness of densities of solutions of stochastic differential equations in
$\mathbb{R}^{n}$ driven by Brownian motion. The domain of definition of the
Malliavin derivative is a subspace $\mathbb{D}_{1,2}$ of $\mathbb{L}^{2}(P)$.
We refer to Nulart \cite{N}, Sanz-Solè \cite{S} and Di Nunno \textit{et al}
\cite{DOP} for information about the Malliavin derivative $D_{t}$ for Brownian
motion and, more generally, Lévy processes. Subsequently, in Aase \textit{et
al} \cite{AaOPU} the Malliavin derivative was put into the context of the
white noise theory of Hida and extended to an operator defined on the whole of
$\mathbb{L}^{2}(P)$ and with values in the Hida space $(\mathcal{S})^{\ast}$
of stochastic distributions. This extension is called the
\emph{Hida-Malliavin} derivative. \newline There are several advantages with
working with this extended Hida-Malliavin derivative:

\begin{itemize}
\item The Hida-Malliavin derivative is defined on all of $\mathbb{L}^{2}(P)$,
and it is an extension of the classical Malliavin derivative, in the sense
that it coincides with the classical Malliavin derivative on the subspace
$\mathbb{D}_{1,2}.$

\item The Hida-Malliavin derivative combines well with the white noise
calculus, including the Skorohod integral and calculus with the Wick product
$\diamond$.
\end{itemize}

Moreover, it extends easily to a Hida-Malliavin derivative with respect to a
Poisson random measure.\newline In the following, we let $(\mathcal{S})^{\ast
}$ denote the Hida space of stochastic distributions.\newline\newline It was
proved in Aase \textit{et al} \cite{AaOPU} that one can extend the
Hida-Malliavin derivative operator $D_{t}$ from $\mathbb{D}_{1,2}$ to all of
$\mathbb{L}^{2}(\mathcal{F}_{T},P)$ in such a way that, also denoting the
extended operator by $D_{t}$, for all $F\in\mathbb{L}^{2}(\mathcal{F}_{T},P)$,
we have
\begin{equation}
D_{t}F\in(\mathcal{S})^{\ast}\text{ and }(t,\omega)\mapsto\mathbb{E}%
[D_{t}F\mid\mathcal{F}_{t}]\text{ belongs to }\mathbb{L}^{2}(\lambda\times P),
\label{eq2.10a}%
\end{equation}
where $\lambda$ is Lebesgue measure on $[0,T].$

\begin{proposition}
[Generalized Clark-Ocone formula \cite{DOP}]For all $F\in\mathbb{L}%
^{2}(\mathcal{F}_{T},P)$, we have%
\begin{equation}
F=\mathbb{E}[F]+%
{\textstyle\int_{0}^{T}}
\mathbb{E}[D_{t}F\mid\mathcal{F}_{t}]dB(t). \label{eq2.11a}%
\end{equation}

\end{proposition}

Moreover, we have the following \emph{generalized duality formula,} for the
Brownian motion:

\begin{proposition}
[The generalized duality formula for $B$]Fix $s\in\lbrack0,T]$. If
$t\mapsto\varphi(t,s,\omega)\in\mathbb{L}^{2}(\lambda\times P)$ is
$\mathbb{F}$-adapted with $\mathbb{E}[%
{\textstyle\int_{0}^{T}}
\varphi^{2}(t,s)dt]<\infty$ and $F\in\mathbb{L}^{2}(\mathcal{F}_{T},P)$, then
we have,
\begin{equation}
\mathbb{E}[F%
{\textstyle\int_{0}^{T}}
\varphi(t,s)dB(t)]=\mathbb{E}[%
{\textstyle\int_{0}^{T}}
\mathbb{E}[D_{t}F\mid\mathcal{F}_{t}]\varphi(t,s)dt]. \label{geduB}%
\end{equation}

\end{proposition}

\noindent{Proof.} \quad As it has observed by Agram and Øksendal \cite{AO},
for fixed $s\in\lbrack0,T],$ by \eqref{eq2.10a}-\eqref{eq2.11a} and the Itô
isometry, we get
\begin{align*}
&  \mathbb{E}[F%
{\textstyle\int_{0}^{T}}
\varphi(t,s)dB(t)]=\mathbb{E}[(\mathbb{E}[F]+%
{\textstyle\int_{0}^{T}}
\mathbb{E}[D_{t}F\mid\mathcal{F}_{t}]dB(t))(%
{\textstyle\int_{0}^{T}}
\varphi(t,s)dB(t))]\\
&  =\mathbb{E}[%
{\textstyle\int_{0}^{T}}
\mathbb{E}[D_{t}F\mid\mathcal{F}_{t}]\varphi(t,s)dt].
\end{align*}
\hfill$\square$ \newline As we have mentioned earlier, there is also an
extension of the Hida-Malliavin derivative with respect to the Poisson random
measure $D_{t,\zeta}$ from $\mathbb{D}_{1,2}^{(\tilde{N})}$ to
$\mathbb{L}^{2}(\mathcal{F}_T, P)$
\color{black}
such that, also denoting the extended operator by $D_{t,\zeta}$, for all
$F\in\mathbb{L}^{2}(\mathcal{F}_{T},P)$, we have
\[
D_{t,\zeta}F\in(\mathcal{S})^{\ast}\text{ and }(t,\zeta,\omega)\mapsto
\mathbb{E}[D_{t,\zeta}F\mid\mathcal{F}_{t}]\text{ belongs to }\mathbb{L}%
^{2}(\lambda\times\nu\times P).
\]
See Di Nunno \textit{et al} \cite{DOP}. Note that in this case, there are two
parameters $t,\zeta,$ where $t\in\lbrack0,T]$ represents time and $\zeta\in%
\mathbb{R}
_{0}$ represents a generic jump size. \newline\newline We now give a jump
diffusion version of the \emph{generalized Clark-Ocone formula}:

\begin{proposition}
[Generalized Clark-Ocone formula \cite{DOP}]For all $F\in%
\mathbb{L}^{2}(\mathcal{F}_{T},P)%
\color{black}%
$, we have%
\[
F=\mathbb{E}[F]+%
{\textstyle\int_{0}^{T}}
{\textstyle\int_{\mathbb{R}_{0}}}
\mathbb{E}[D_{t,\zeta}F|\mathcal{F}_{t}]\tilde{N}(dt,d\zeta),
\]
where we have chosen a predictable version of the conditional expectation
$\mathbb{E}[D_{t,\zeta}F|\mathcal{F}_{t}]$ for each $t\geq0.$
\end{proposition}

Moreover, we have also an extension of the duality formula for jumps:

\begin{proposition}
[The generalized duality formula for $\tilde{N}$]Fix $s\in\lbrack0,T]$.
Suppose $\Psi(t,s,\zeta)$ is $\mathbb{F}$-adapted and $\mathbb{E}[%
{\textstyle\int_{0}^{T}}
{\textstyle\int_{\mathbb{R}_{0}}}
\Psi^{2}(t,s,\zeta)\nu(d\zeta)dt]<\infty$ and let $F\in%
\mathbb{L}^{2}(\mathcal{F}_{T},P)%
\color{black}%
$. Then,
\begin{equation}
\mathbb{E}[F%
{\textstyle\int_{0}^{T}}
{\textstyle\int_{\mathbb{R}_{0}}}
\Psi(t,s,\zeta)\tilde{N}(dt,d\zeta)]=\mathbb{E}[%
{\textstyle\int_{0}^{T}}
{\textstyle\int_{\mathbb{R}_{0}}}
\mathbb{E}[D_{t,\zeta}F|\mathcal{F}_{t}]\Psi(t,s,\zeta)\nu(d\zeta)dt].
\label{geduN}%
\end{equation}
Accordingly, note that from now on we are working with this generalized
version of the Malliavin derivative. We emphasize that this generalized
Hida-Malliavin derivative $DX$ (where $D$ stands for $D_{t}$ or $D_{t,\zeta}$,
depending on the setting) exists for all $X\in\mathbb{L}^{2}(P)$ as an element
of the Hida stochastic distribution space $(\mathcal{S})^{\ast}$, and it has
the property that the conditional expectation $\mathbb{E}[DX|\mathcal{F}_{t}]$
belongs to $\mathbb{L}^{2}(\lambda\times P)$, where $\lambda$ is a Lebesgue
measure on $[0,T]$. Therefore, when using the Hida-Malliavin derivative,
combined with conditional expectation, no assumptions on Hida-Malliavin
differentiability in the classical sense are needed; we can work on the whole
space of random variables in $\mathbb{L}^{2}(P)$.\newline
\end{proposition}

The following result is the Hida-Malliavin representation for BSVIE:

\begin{theorem}
[Representation theorem for BSVIE]\label{Th2.9} Suppose that the driver
$f(t,s,p,q,r):\left[  0,T\right]  ^{2}\times\mathbb{R}\times\mathbb{R\times
L}_{\nu}^{2}\times\Omega\rightarrow\mathbb{R}$ \ is $\mathbb{F}$-adapted with
respect to $s$ for all $t,p,q,r$ and that $s\mapsto(p(s),q(t,s),r(t,s,\zeta
))\in\mathbb{L}^{2}\times\mathbb{L}^{2}\times\mathbb{H}_{\nu}^{2}$ are given
$\mathbb{F}$-adapted processes with respect to $s\in\lbrack t,T]$, and they
satisfy
\begin{equation}%
\begin{array}
[c]{cc}%
p(t) & =F(t)+%
{\textstyle\int_{t}^{T}}
f(t,s,p(s),q(t,s),r(t,s,\cdot))ds-%
{\textstyle\int_{t}^{T}}
q(t,s)dB(s)\\
& -%
{\textstyle\int_{t}^{T}}
{\textstyle\int_{\mathbb{R}_{0}}}
r(t,s,\zeta)\tilde{N}(ds,d\zeta);t\in\lbrack0,T],
\end{array}
\label{BSVIE_J}%
\end{equation}

where $F(t)\in\mathbb{L}^{2}(\mathcal{F}_{T},P,%
\mathbb{R}
).$ Then for a.a. $t,s$ and $\zeta,$ the following holds:
\begin{equation}
q(t,s)=\mathbb{E}[D_{s}p(t)|\mathcal{F}_{s}]; \quad s<t \label{eq2.25}%
\end{equation}

and
\begin{equation}
r(t,s,\zeta)=\mathbb{E}[D_{s,\zeta}p(t)|\mathcal{F}_{s}];\quad s<t.
\label{eq2.27}%
\end{equation}

\end{theorem}

\noindent{Proof.} \quad We know by Theorem 3.1 in Agram \textit{el al}
\cite{AOY} that for a Lipschitz driver $f$ and for a terminal value
$F(t)\in\mathbb{L}^{2}(\mathcal{F}_{T},P,%
\mathbb{R}
),$ the above BSVIE with jumps (\ref{BSVIE_J}) has a unique solution.
Moreover, for all $t\in\lbrack0,T],$ it holds that%
\begin{equation}
p(t)=\mathbb{E}[p(t)]+%
{\textstyle\int_{0}^{t}}
q(t,s)dB(s)+%
{\textstyle\int_{0}^{t}}
{\textstyle\int_{\mathbb{R}_{0}}}
r(t,s,\zeta)\tilde{N}(ds,d\zeta);t\in\lbrack0,T]. \label{est}%
\end{equation}
For more details, we refer to Yong \cite{Yong1} (for the Brownian framework)
and for the discountinuous case, we refer to Ren \cite{Ren}.\newline Taking
the Hida-Malliavin derivatives of \eqref{est}, we get, for $s<t,$
\begin{equation}
\mathbb{E}[D_{s}p(t)|\mathcal{F}_{s}]=q(t,s)+\mathbb{E}[{\textstyle\int
_{s}^{t}}D_{s}p(t)dB(s)|\mathcal{F}_{s}]=q(t,s).
\end{equation}
Similarly, for $s<t$,
\[
\mathbb{E}[D_{s,\zeta}p(t)|\mathcal{F}_{s}]=r(t,s,\zeta)+\mathbb{E}%
[{\textstyle\int_{s}^{t}}{\textstyle\int_{\mathbb{R}_{0}}}D_{s,\zeta
}r(t,s,\zeta))\tilde{N}(ds,d\zeta)|\mathcal{F}_{s}]=r(t,s,\zeta).
\]

\hfill$\square$ \bigskip

\section{Existence and uniqueness of solutions of SVIEs}

In this section, we prove existence and uniqueness of solutions of SVIE driven
by Brownian motion and an independent compensated Poisson random measure,
under some conditions. The case of SVIE driven by right continuous
semimartingales in general has been studied by Protter \cite{p}.\newline Let
us consider the stochastic Volterra integral equation with jumps of the form%

\begin{equation}%
\begin{array}
[c]{l}%
X(t)=X^{u}(t)=\xi(t)+%
{\textstyle\int_{0}^{t}}
b\left(  t,s,X(s),u(s)\right)  ds+%
{\textstyle\int_{0}^{t}}
\sigma\left(  t,s,X(s),u(s)\right)  dB(s)\\
\text{ \ \ \ \ \ \ \ \ \ \ \ }+%
{\textstyle\int_{0}^{t}}
{\textstyle\int_{\mathbb{R}_{0}}}
\gamma\left(  t,s,X(s),u(s),\zeta\right)  \tilde{N}(ds,d\zeta);t\in
\lbrack0,T].
\end{array}
\label{SVIE}%
\end{equation}
We impose the following set of assumptions:

\begin{description}
\item[(i)] $\xi(t)$ is a given $\mathbb{F}$-adapted càdlàg process,

\item[(ii)] $b\left(  t,s,x,u\right)  ,$ $\sigma(t,s,x,u)$ and $\gamma
(t,s,x,u,\zeta)$ are $\mathcal{F}_{s}$-predictable for $s\leq t$, for each
$(t,x,u)$ and $(t,x,u,\zeta)$ respectively.

\item[(iii)] There exists a constant $C>0,$ such that, for all $0\leq s\leq
t\leq T,\,u\in U$ and all $b,\sigma,\gamma$ satisfy
\end{description}

\[
\left\vert b(t,s,0,u)\right\vert +\left\vert \sigma(t,s,0,u)\right\vert +(%
{\textstyle\int_{\mathbb{R}_{0}}}
\left\vert \gamma(t,s,0,u,\zeta)\right\vert ^{2}\nu(d\zeta))^{\frac{1}{2}}\leq
C,
\]

\begin{description}
\item[(iv)] $b(t,s,\cdot,u),\sigma(t,s,\cdot,u)$ and $\gamma(t,s,\cdot
,u,\zeta)$ are Lipschitz continuous with respect to $x$ uniformly in $t,s,u,$
i.e. for all $x,x^{\prime}\in%
\mathbb{R}
$, we have
\[%
\begin{array}
[c]{c}%
\left\vert b(t,s,x,u)-b(t,s,x^{\prime},u)\right\vert +\left\vert
\sigma(t,s,x,u)-\sigma(t,s,x^{\prime},u)\right\vert \\
+(%
{\textstyle\int_{\mathbb{R}_{0}}}
\left\vert \gamma(t,s,x,u,\zeta)-\gamma(t,s,x^{\prime},u,\zeta)\right\vert
^{2}\nu(d\zeta))^{\frac{1}{2}}\\
\leq C|x-x^{\prime}|,
\end{array}
\]

\item[(v)]
\[%
\begin{array}
[c]{c}%
\left\vert b(t,s,x,u)\right\vert +\left\vert \sigma(t,s,x,u)\right\vert +(%
{\textstyle\int_{\mathbb{R}_{0}}}
\left\vert \gamma(t,s,x,u,\zeta)\right\vert ^{2}\nu(d\zeta))^{\frac{1}{2}}\\
\leq C(|1+|x|).
\end{array}
\]

\end{description}

\begin{theorem}
Under the above assumptions $\left(  i-v\right)  $, the SVIE (\ref{SVIE}) has
a unique solution.
\end{theorem}

\noindent{Proof.} \quad\textbf{Existence}. Fix $u\in U.$ Define $X^{n}$
inductively for $n=0,1$,.., as follows%
\begin{equation}
\left\{
\begin{array}
[c]{ll}%
X^{0}(t) & =\xi(t),\\
X^{n+1}(t) & =\xi(t)+%
{\textstyle\int_{0}^{t}}
b\left(  t,s,X^{n}(s),u(s)\right)  ds+%
{\textstyle\int_{0}^{t}}
\sigma\left(  t,s,X^{n}(s),u(s)\right)  dB(s)\\
& +%
{\textstyle\int_{0}^{t}}
{\textstyle\int_{\mathbb{R}_{0}}}
\gamma\left(  t,s,X^{n}(s),u(s),\zeta\right)  \tilde{N}(ds,d\zeta);t\in\left[
0,T\right]  ,\,n\geq0.
\end{array}
\right.  \label{ex}%
\end{equation}
Let $\overline{X}^{n}:=X^{n+1}-X^{n}.$ Then, the following estimate holds, for
each $t\in\left[  0,T\right]  ,$ and $n\geq1,$%
\begin{equation}%
\begin{array}
[c]{l}%
\mathbb{E}[\left\vert \overline{X}^{n}(s)\right\vert ^{2}]\leq3\mathbb{E}[t%
{\textstyle\int_{0}^{t}}
|b\left(  t,s,X^{n}(s),u(s)\right)  -b\left(  t,s,X^{n-1}(s),u(s)\right)
|^{2}ds\\
+%
{\textstyle\int_{0}^{t}}
|\sigma\left(  t,s,X^{n}(s),u(s)\right)  -\sigma\left(  t,s,X^{n-1}%
(s),u(s)\right)  |^{2}ds\\
+%
{\textstyle\int_{0}^{t}}
{\textstyle\int_{\mathbb{R}_{0}}}
|\gamma\left(  t,s,X^{n}(s),u(s),\zeta\right)  -\gamma\left(  t,s,X^{n-1}%
(s),u(s),\zeta\right)  |^{2}\nu(d\zeta)ds].
\end{array}
\label{3}%
\end{equation}
Using assumption (iv), we get%
\[
\mathbb{E}[\left\vert \overline{X}^{n}(t)\right\vert ^{2}]\leq3C^{2}%
\mathbb{E}[t%
{\textstyle\int_{0}^{t}}
\{|\overline{X}^{n-1}(s)|^{2}+2|\overline{X}^{n-1}(s)|^{2}\}ds].
\]
Define $K:=3C^{2}(T+2)$, then%
\begin{equation}
\mathbb{E}[\left\vert X^{n+1}(t)-X^{n}(t)\right\vert ^{2}]\leq K\mathbb{E}[%
{\textstyle\int_{0}^{t}}
|X^{n}(s)-X^{n-1}(s)|^{2}ds]. \label{xn}%
\end{equation}
Now, using the linear growth assumption (v), we obtain similarly as above%
\[%
\begin{array}
[c]{l}%
\mathbb{E}[\left\vert X^{1}(t)-X^{0}(t)\right\vert ^{2}]=\mathbb{E}[|%
{\textstyle\int_{0}^{t}}
b\left(  t,s,\xi(s),u(s)\right)  ds\\
+%
{\textstyle\int_{0}^{t}}
\sigma\left(  t,s,\xi(s),u(s)\right)  dB(s)\\
+%
{\textstyle\int_{0}^{t}}
{\textstyle\int_{\mathbb{R}_{0}}}
\gamma\left(  t,s,\xi(s),u(s),\zeta\right)  \tilde{N}(ds,d\zeta)|^{2}]\\
\leq2K\mathbb{E}%
{\textstyle\int_{0}^{t}}
[1+\xi^{2}(s)]ds\\
\leq2Kt\underset{t\in\left[  0,T\right]  }{\sup}\mathbb{E}[1+\xi^{2}(t)].
\end{array}
\]
Combine this with (\ref{xn}), yields%
\[
\mathbb{E}[\left\vert X^{n+1}(t)-X^{n}(t)\right\vert ^{2}]\leq\tfrac
{2K^{\prime}(Kt)^{n+1}}{(n+1)!},
\]
where $K^{\prime}:=\underset{t\in\left[  0,T\right]  }{\sup}\mathbb{E}%
[1+\xi^{2}(t)]<\infty.$ For $m>n>0,$ it follows that%

\begin{align*}
\mathbb{E}[%
{\textstyle\int_{0}^{T}}
\left\vert X^{m}(t)-X^{n}(t)\right\vert ^{2}dt]  &  \leq\overset
{m-1}{\underset{k=n}{%
{\textstyle\sum}
}}\tfrac{2K^{\prime}%
{\textstyle\int_{0}^{T}}
(Kt)^{k+1}dt}{(k+1)!}\\
&  =\overset{m-1}{\underset{k=n}{%
{\textstyle\sum}
}}\tfrac{2K^{\prime}K^{k+1}T^{k+2}}{(k+2)!}\rightarrow0,\text{ as
}m,n\rightarrow\infty.
\end{align*}
Hence, $\{X^{n}(t)\}_{n=1}^{\infty}$ is a Cauchy sequence in $\mathbb{L}%
^{2}(\lambda\times P)$. Finally, taking the limit in the Picard iteration as
$n\rightarrow+\infty,$ yields
\[%
\begin{array}
[c]{l}%
X(t)=\xi(t)+%
{\textstyle\int_{0}^{t}}
b\left(  t,s,X(s),u(s)\right)  ds+%
{\textstyle\int_{0}^{t}}
\sigma\left(  t,s,X(s),u(s)\right)  dB(s)\\
\text{ \ \ \ \ \ \ \ \ \ \ \ }+%
{\textstyle\int_{0}^{t}}
{\textstyle\int_{\mathbb{R}_{0}}}
\gamma\left(  t,s,X(s),u(s),\zeta\right)  \tilde{N}(ds,d\zeta);t\in
\lbrack0,T].
\end{array}
\]
\textbf{Uniqueness.} The uniqueness is obtained by the estimate of the
difference of two solutions, and it is carried out similarly to the argument
above.$\qquad\qquad\qquad\qquad\qquad\qquad\qquad\qquad\qquad\qquad\square$

\section{Stochastic maximum principles}

In this section, we study stochastic maximum principles of stochastic Volterra
integral systems under partial information, i.e., the information available to
the controller is given by a sub-filtration $\mathbb{G}=\{\mathcal{G}%
_{t}\}_{t\geq0}$ such that $\mathcal{G}_{t}\subseteq\mathcal{F}_{t}$ for all
$t\geq0.$ The set $U\subset\mathbb{R}$ is assumed to be convex. The set of
admissible controls, i.e. the strategies available to the controller is given
by a subset $\mathcal{A}_{\mathbb{G}}$ of the càdlàg, $U$-valued and
$\mathbb{G}$-adapted processes.\newline The state of our system $X^{u}%
(t)=X(t)$ satisfies the following SVIE%
\begin{equation}%
\begin{array}
[c]{cc}%
X(t) & =\xi(t)+%
{\textstyle\int_{0}^{t}}
b\left(  t,s,X(s),u(s)\right)  ds+%
{\textstyle\int_{0}^{t}}
\sigma\left(  t,s,X(s),u(s)\right)  dB(s)\\
& +%
{\textstyle\int_{0}^{t}}
{\textstyle\int_{\mathbb{R}_{0}}}
\gamma\left(  t,s,X(s),u(s),\zeta\right)  \tilde{N}(ds,d\zeta);t\in
\lbrack0,T],
\end{array}
\label{sde}%
\end{equation}
where $b(t,s,x,u)=b(t,s,x,u,\omega):\left[  0,T\right]  ^{2}\times
\mathbb{R}\times U\times\Omega\rightarrow\mathbb{R}$, $\sigma(t,s,x,u)=\sigma
(t,s,x,u,\omega):\left[  0,T\right]  ^{2}\times\mathbb{R}\times U\times
\Omega\rightarrow%
\mathbb{R}
$ and $\gamma(t,s,x,u,\zeta)=\gamma(t,s,x,u,\zeta,\omega):\left[  0,T\right]
^{2}\times\mathbb{R}\times U\times%
\mathbb{R}
_{0}\times\Omega\rightarrow%
\mathbb{R}
$.\newline\newline The \emph{performance functional} has the form%
\begin{equation}%
\begin{array}
[c]{ll}%
J(u) & =\mathbb{E[}{%
{\textstyle\int_{0}^{T}}
}\text{ }f(t,X(t),u(t))dt+g(X(T))],\quad u\in\mathcal{A}_{\mathbb{G}},
\end{array}
\label{perf}%
\end{equation}
with given functions $f(t,x,u)=f(t,x,u,\omega):\left[  0,T\right]
\times\mathbb{R}\times U\times\Omega\rightarrow\mathbb{R}$ and
$g(x)=g(x,\omega):\mathbb{R\times}\Omega\rightarrow\mathbb{R}$.\\

 We
impose the following assumption:\\

\emph{Assumption A1} \newline \emph{The
processes $b,\sigma,f$ and $\gamma$ are $\mathcal{F}_{s}$-adapted for all
$s\leq t,$ and twice continuously differentiable ($C^{2}$) with respect to
$t$, $x$ and continuously differentiable ($C^{1}$) with respect to $u$ for
each $s$. The driver $g$ is assumed to be $\mathcal{F}_{T}$-measurable and
$C^{1}$ in $x$. Moreover, all the partial derivatives are supposed to be
bounded.\\
}

 Note that the performance functional (\ref{perf}) is not of
Volterra type.

\subsection{The Hamiltonian and the adjoint equations}

Define the \emph{Hamiltonian functional} associated to our control problem
(\ref{sde}) and (\ref{perf}), as
\begin{equation}%
\begin{array}
[c]{l}%
\mathcal{H}(t,x,v,p,q,r(\cdot))\\
:=H^{0}(t,x,v,p,q,r(\cdot))+H^{1}(t,x,v,p,q,r(\cdot)),
\end{array}
\label{eq3.3}%
\end{equation}
where%

\[
H^{0}:[0,T]\times\mathbb{R}\times U\times\mathbb{R}\times\mathbb{R}%
\times\mathbb{L}_{\nu}^{2}\rightarrow\mathbb{R}%
\]
and
\[
H^{1}:[0,T]\times\mathbb{R}\times U\times\mathbb{R}\times\mathbb{R}%
\times\mathbb{L}_{\nu}^{2}\rightarrow\mathbb{R}%
\]
by
\begin{align}
H^{0}(t,x,v,p,q,r(\cdot))  &  :=f(t,x,v)+p(t)b(t,t,x,v)+q(t,t)\sigma
(t,t,x,v)\label{eq2.1}\\
&  +%
{\textstyle\int_{\mathbb{R}_{0}}}
r(t,t,\zeta)\gamma(t,t,x,v,\zeta)\nu(d\zeta),\nonumber
\end{align}%
\begin{align*}
H^{1}(t,x,v,p,q,r(\cdot))  &  :=%
{\textstyle\int_{t}^{T}}
p(s)\tfrac{\partial b}{\partial s}(s,t,x,v)ds+%
{\textstyle\int_{t}^{T}}
q(s,t)\tfrac{\partial\sigma}{\partial s}(s,t,x,v)ds\\
&  +%
{\textstyle\int_{t}^{T}}
{\textstyle\int_{\mathbb{R}_{0}}}
r(s,t,\zeta)\tfrac{\partial\gamma}{\partial s}(s,t,x,v,\zeta)\nu(d\zeta)ds.
\end{align*}
We may regard $x,p,q,r=r(\cdot)$ as generic values for the processes $X(\cdot),$
$p(\cdot),$ $q(\cdot),$ $r(\cdot)$, respectively.\newline The BSVIE for the
adjoint processes $p(t),q(t,s),r(t,s,\cdot)$ is defined by%
\begin{equation}%
\begin{array}
[c]{cc}%
p(t) & =\tfrac{\partial g}{\partial x}(X(T))+%
{\textstyle\int_{t}^{T}}
\tfrac{\partial\mathcal{H}}{\partial x}(s)ds-%
{\textstyle\int_{t}^{T}}
q(t,s)dB(s)\\
& -%
{\textstyle\int_{t}^{T}}
{\textstyle\int_{\mathbb{R}_{0}}}
r(t,s,\zeta)\tilde{N}(ds,d\zeta);t\in\lbrack0,T],
\end{array}
\label{p}%
\end{equation}
where we have used the simplified notation
\[
\tfrac{\partial\mathcal{H}}{\partial x}(t)=\tfrac{\partial\mathcal{H}%
}{\partial x}(t,X(t),u(t),p(t),q(t,t),r(t,t,\cdot)).
\]

\begin{remark}
Using the definition of $\mathcal{H}$ and the Fubini theorem, we see that the
driver in the BSVIE \eqref{p} can be explicitly written
\begin{align}
&  \int_{t}^{T}\frac{\partial\mathcal{H}}{\partial x}(s)ds=\int_{t}%
^{T}\Big\{\frac{\partial f}{\partial x}(s,x,v)+p(s)\frac{\partial b}{\partial
x}(s,s,x,v)+\int_{s}^{T}p(z)\frac{\partial^{2}b}{\partial z\partial
x}(z,t,x,v)dz\nonumber\\
&  +q(s,s)\frac{\partial\sigma}{\partial x}(s,s,x,v)+\int_{s}^{T}%
q(z,t)\frac{\partial^{2}\sigma}{\partial z\partial x}(z,t,x,v)dz\nonumber\\
&  +\int_{\mathbb{R}_{0}}r(s,s,\zeta)\frac{\partial\gamma}{\partial
x}(s,s,x,v,\zeta)\nu(d\zeta)+\int_{s}^{T}\int_{\mathbb{R}_{0}}r(z,t,\zeta
)\frac{\partial^{2}\gamma}{\partial z\partial x}(z,t,x,v,\zeta)\nu
(d\zeta)dz\Big\}ds\nonumber\\
&  =\int_{t}^{T}\Big\{\frac{\partial f}{\partial x}(s,x,v)+p(s)\Big[\frac
{\partial b}{\partial x}(s,s,x,v)+
(s-t)
\color{black}%
\frac{\partial^{2}b}{\partial s\partial x}(s,t,x,v)\Big]\nonumber\\
&  +q(s,s)\frac{\partial\sigma}{\partial x}(s,s,x,v)+%
(s-t)%
\color{black}%
q(s,t)\frac{\partial^{2}\sigma}{\partial s\partial x}(s,t,x,v)\nonumber\\
&  +\int_{\mathbb{R}_{0}}\Big[r(s,s,\zeta)\frac{\partial\gamma}{\partial
x}(s,s,x,v,\zeta)+
(s-t)
\color{black}%
r(s,t,x,v,\zeta)\frac{\partial^{2}\gamma}{\partial s\partial x}(s,t,x,v,\zeta
)\Big]\nu(d\zeta)\Big\}ds.
\end{align}

\end{remark}

From this it follows by Theorem 3.1 in Agram \textit{et al }\cite{AOY}, that
we have existence and uniqueness of the solution of equation \eqref{p}.
\newline%
From now on we also make the following assumption:\\

\emph{Assumption }$\emph{A2}$\newline \emph{The functions $t\mapsto
q(t,s)$ and $t\mapsto r(t,s,\cdot)$ are $C^{1}$ for all $s,\zeta,\omega$ and
\[
\mathbb{E}[%
{\textstyle\int_{0}^{T}}
{\textstyle\int_{0}^{T}}
(\tfrac{\partial q(t,s)}{\partial t})^{2}dsdt+%
{\textstyle\int_{0}^{T}}
{\textstyle\int_{0}^{T}}
{\textstyle\int_{\mathbb{R}_{0}}}
(\tfrac{\partial r(t,s,\zeta)}{\partial t})^{2}\nu(d\zeta)dsdt]<\infty.
\]
}
\vskip 0.3cm

Note that from equation (\ref{sde}), we get the following equivalent
formulation, for each $(t,s)\in\lbrack0,T]^{2},$
\begin{equation}%
\begin{array}
[c]{ll}%
dX(t) & =\xi^{\prime}(t)dt+b\left(  t,t,X(t),u(t)\right)  dt+(%
{\textstyle\int_{0}^{t}}
\tfrac{\partial b}{\partial t}\left(  t,s,X(s),u(s)\right)  ds)dt\\
& +\sigma\left(  t,t,X(t),u(t)\right)  dB(t)+(%
{\textstyle\int_{0}^{t}}
\tfrac{\partial\sigma}{\partial t}\left(  t,s,X(s),u(s)\right)  dB(s))dt\\
& +%
{\textstyle\int_{\mathbb{R}_{0}}}
\gamma\left(  t,t,X(t),u(t),\zeta\right)  \tilde{N}(dt,d\zeta)+(%
{\textstyle\int_{0}^{t}}
{\textstyle\int_{\mathbb{R}_{0}}}
\tfrac{\partial\gamma}{\partial t}\left(  t,s,X(s),u(s),\zeta\right)
\tilde{N}(ds,d\zeta))dt,
\end{array}
\label{eq3.6}%
\end{equation}
and from equation (\ref{p}) under assumption $\emph{A2}$, 
we have the following differential form
\begin{equation}
\left\{
\begin{array}
[c]{ll}%
dp(t) & =-[\tfrac{\partial\mathcal{H}}{\partial x}(t)+%
{\textstyle\int_{t}^{T}}
\tfrac{\partial q}{\partial t}(t,s)dB(s)+%
{\textstyle\int_{t}^{T}}
{\textstyle\int_{\mathbb{R}_{0}}}
\tfrac{\partial r}{\partial t}(t,s,\zeta)\tilde{N}(ds,d\zeta)]dt\\
& +q(t,t)dB(t)+%
{\textstyle\int_{\mathbb{R}_{0}}}
r(t,t,\zeta)\tilde{N}(dt,d\zeta),\\
p(T) & =\tfrac{\partial g}{\partial x}(X(T)).
\end{array}
\right.  \label{eq3.7}%
\end{equation}

\begin{remark}
Assumption $\emph{A2}$ is verified in a subclass of linear BSVIE with
jumps, as we will see in section 5. For more details, we refer to Hu and
Øksendal \cite{ho}.
\end{remark}

\subsection{A sufficient maximum principle}

We now state and prove a sufficient version of the maximum principle approach
(a verification theorem).

\begin{theorem}
[Sufficient maximum principle]Let $\hat{u}\in\mathcal{A}_{\mathbb{G}},$ with
corresponding solutions $\hat{X}(t),$ $\left(  \hat{p}(t),\hat{q}(t,s),\hat
{r}(t,s,\cdot)\right)  $ of (\ref{sde}) and (\ref{p}) respectively. Assume that

\begin{itemize}
\item The functions%
\[
x\mapsto g(x),
\]
and
\[
(x,u)\mapsto\mathcal{H}(t,x,u,\hat{p},\hat{q},\hat{r}(\cdot))
\]
are concave.

\item (The maximum condition)
\begin{align}
&  \underset{v\in\mathbb{U}}{\sup}\text{ }\mathbb{E[}\mathcal{H}(t,\hat
{X}(t),v,\hat{p}(t),\hat{q}(t,t),\hat{r}(t,t,\cdot))|{\mathcal{G}}%
_{t}]\nonumber\\
&  =\mathbb{E[}\mathcal{H}(t,\hat{X}(t),\hat{u}(t),\hat{p}(t),\hat
{q}(t,t),\hat{r}(t,t,\cdot))|{\mathcal{G}}_{t}]\text{ }\forall t\text{
}P\text{-a.s.} \label{eq3.9}%
\end{align}
\newline Then $\hat{u}$ is an optimal control for our problem.
\end{itemize}
\end{theorem}

\noindent{Proof.} \quad By considering a sequence of stopping times converging
upwards to $T$, we see that we may assume that all the $dB$- and $\tilde{N}$-
integrals in the following are martingales and hence have expectation $0$. We
refer to the proof of Lemma 3.1 in \cite{OS} for details. \newline\newline
Choose $u\in\mathcal{A}_{\mathbb{G}},$ we want to prove that $J(u)\leq
J(\hat{u})$.\newline By the definition of the cost functional (\ref{perf})$,$
we have
\begin{equation}
J(u)-J(\hat{u})=I_{1}+I_{2}, \label{eq2.8}%
\end{equation}
where we have used the shorthand notations \newline%
\[
I_{1}=\mathbb{E[}%
{\textstyle\int_{0}^{T}}
\tilde{f}\left(  t\right)  dt],\quad I_{2}=\mathbb{E[}\tilde{g}(T)],
\]
\newline and
\[
\tilde{f}\left(  t\right)  =f(t)-\hat{f}(t),
\]
with
\[%
\begin{array}
[c]{ll}%
f(t) & =f\left(  t,X(t),u(t)\right)  ,\\
\hat{f}\left(  t\right)  & =f(t,\hat{X}(t),\hat{u}(t)),
\end{array}
\]
and similarly for $b(t,t)=b\left(  t,t,X(t),u(t)\right)  ,$ and the other
coefficients. By the definition of the Hamiltonian $(\ref{eq2.1})$, we get
\begin{equation}
I_{1}=\mathbb{E}[%
{\textstyle\int_{0}^{T}}
\{\tilde{H}^{0}(t)-\hat{p}(t)\tilde{b}(t,t)-\hat{q}(t,t)\tilde{\sigma}(t,t)-%
{\textstyle\int_{\mathbb{R}_{0}}}
\hat{r}(t,t,\zeta)\tilde{\gamma}(t,t,\zeta)\nu(d\zeta)\}dt], \label{i1}%
\end{equation}
where $\tilde{H}^{0}(t)=H^{0}(t)-\hat{H}^{0}(t)$ with%
\[%
\begin{array}
[c]{ll}%
H^{0}(t) & =H^{0}(t,X(t),u(t),\hat{p}(t),\hat{q}(t,t),\hat{r}(t,t,\cdot)),\\
\hat{H}^{0}(t) & =\hat{H}^{0}(t,\hat{X}(t),\hat{u}(t),\hat{p}(t),\hat
{q}(t,t),\hat{r}(t,t,\cdot)).
\end{array}
\]
By the concavity of $g$ and the terminal value of the BSVIE $\left(
\ref{p}\right)  $, we obtain%
\[%
\begin{array}
[c]{lll}%
I_{2} & \leq\mathbb{E}[\tfrac{\partial\hat{g}}{\partial x}(T)\tilde{X}(T)] &
=\mathbb{E}[\hat{p}(T)\tilde{X}(T)].
\end{array}
\]
Applying the Itô formula to $\hat{p}(t)\tilde{X}(t)$, we get
\begin{align}
I_{2}  &  \leq\mathbb{E}[\hat{p}(T)\tilde{X}(T)]\nonumber\\
&  =\mathbb{E[}%
{\textstyle\int_{0}^{T}}
\hat{p}(t)\{\tilde{b}(t,t)+%
{\textstyle\int_{0}^{t}}
\tfrac{\partial\tilde{b}}{\partial t}(t,s)ds+%
{\textstyle\int_{0}^{t}}
\tfrac{\partial\tilde{\sigma}}{\partial t}\left(  t,s\right)  dB(s)\nonumber\\
&  +%
{\textstyle\int_{0}^{t}}
{\textstyle\int_{\mathbb{R}_{0}}}
\tfrac{\partial\tilde{\gamma}}{\partial t}\left(  t,s,\zeta\right)  \tilde
{N}(ds,d\zeta)\}dt+%
{\textstyle\int_{0}^{T}}
\tilde{X}(t)\{-\tfrac{\partial\widehat{\mathcal{H}}}{\partial x}(t)+%
{\textstyle\int_{t}^{T}}
\tfrac{\partial\hat{q}}{\partial t}(t,s)dB(s)\nonumber\\
&  +%
{\textstyle\int_{t}^{T}}
{\textstyle\int_{\mathbb{R}_{0}}}
\tfrac{\partial\hat{r}}{\partial t}(t,s,\zeta)\tilde{N}(ds,d\zeta)\}dt+%
{\textstyle\int_{0}^{T}}
\hat{q}(t,t)\tilde{\sigma}(t,t)dt+%
{\textstyle\int_{0}^{T}}
{\textstyle\int_{\mathbb{R}_{0}}}
\hat{r}(t,t,\zeta)\tilde{\gamma}(t,t,\zeta)\nu(d\zeta)dt]. \label{I2}%
\end{align}
By the Fubini theorem, we get
\begin{equation}%
{\textstyle\int_{0}^{T}}
\hat{p}(t)(%
{\textstyle\int_{0}^{t}}
\tfrac{\partial\tilde{b}}{\partial t}(t,s)ds)dt=%
{\textstyle\int_{0}^{T}}
(%
{\textstyle\int_{s}^{T}}
\hat{p}(t)\tfrac{\partial\tilde{b}}{\partial t}(t,s)dt)ds=%
{\textstyle\int_{0}^{T}}
(%
{\textstyle\int_{t}^{T}}
\hat{p}(s)\tfrac{\partial\tilde{b}}{\partial s}(s,t)ds)dt. \label{eq2.11}%
\end{equation}
The generalized duality formula for the Brownian motion $\left(
\ref{geduB}\right)  $, yields
\begin{align*}
\mathbb{E}[%
{\textstyle\int_{0}^{T}}
\hat{p}(t)(%
{\textstyle\int_{0}^{t}}
\tfrac{\partial\tilde{\sigma}}{\partial t}(t,s)dB(s))dt]  &  =%
{\textstyle\int_{0}^{T}}
\mathbb{E}[%
{\textstyle\int_{0}^{t}}
\hat{p}(t)\tfrac{\partial\tilde{\sigma}}{\partial t}(t,s)dB(s)]dt\\
&  =%
{\textstyle\int_{0}^{T}}
\mathbb{E}[%
{\textstyle\int_{0}^{t}}
\mathbb{E}[D_{s}\hat{p}(t)|\mathcal{F}_{s}]\tfrac{\partial\tilde{\sigma}%
}{\partial t}(t,s)ds]dt.
\end{align*}
Fubini's theorem gives%
\begin{align*}
\mathbb{E}[%
{\textstyle\int_{0}^{T}}
\hat{p}(t)(%
{\textstyle\int_{0}^{t}}
\tfrac{\partial\tilde{\sigma}}{\partial t}(t,s)dB(s))dt]  &  =%
{\textstyle\int_{0}^{T}}
\mathbb{E}[%
{\textstyle\int_{s}^{T}}
\mathbb{E}[D_{s}\hat{p}(t)|\mathcal{F}_{s}]\tfrac{\partial\tilde{\sigma}%
}{\partial t}(t,s)dt]ds\\
&  =\mathbb{E}[%
{\textstyle\int_{0}^{T}}
{\textstyle\int_{t}^{T}}
\mathbb{E}[D_{t}\hat{p}(s)|\mathcal{F}_{t}]\tfrac{\partial\tilde{\sigma}%
}{\partial s}(s,t)dsdt],
\end{align*}
and by equality $\left(  \ref{eq2.25}\right)  $, we end up with%
\begin{equation}
\mathbb{E}[%
{\textstyle\int_{0}^{T}}
\hat{p}(t)(%
{\textstyle\int_{0}^{t}}
\tfrac{\partial\tilde{\sigma}}{\partial t}(t,s)dB(s))dt]=\mathbb{E}[%
{\textstyle\int_{0}^{T}}
{\textstyle\int_{t}^{T}}
\hat{q}(s,t)\tfrac{\partial\tilde{\sigma}}{\partial s}(s,t)dsdt].
\label{eq2.12}%
\end{equation}
Doing similar considerations as for the Brownian setting for the jumps, such
as the Fubini theorem, the generalized duality formula for jumps (\ref{geduN})
and $\left(  \ref{eq2.27}\right)  $, we obtain%

\begin{align}
\mathbb{E}[%
{\textstyle\int_{0}^{T}}
(%
{\textstyle\int_{0}^{t}}
{\textstyle\int_{\mathbb{R}_{0}}}
\hat{p}(t)\tfrac{\partial\tilde{\gamma}}{\partial t}(t,s,\zeta)\tilde
{N}(ds,d\zeta))dt]  &  =%
{\textstyle\int_{0}^{T}}
\mathbb{E}[%
{\textstyle\int_{0}^{t}}
{\textstyle\int_{\mathbb{R}_{0}}}
\hat{p}(t)\tfrac{\partial\tilde{\gamma}}{\partial t}(t,s,\zeta)\tilde
{N}(ds,d\zeta))]dt\nonumber\\
&  =%
{\textstyle\int_{0}^{T}}
\mathbb{E}[%
{\textstyle\int_{0}^{t}}
{\textstyle\int_{\mathbb{R}_{0}}}
\mathbb{E}[D_{s,\zeta}\hat{p}(t)|\mathcal{F}_{s}]\tfrac{\partial\tilde{\gamma
}}{\partial t}(t,s,\zeta)\nu(d\zeta)ds]dt\nonumber\\
&  =%
{\textstyle\int_{0}^{T}}
\mathbb{E}[%
{\textstyle\int_{s}^{T}}
{\textstyle\int_{\mathbb{R}_{0}}}
\mathbb{E}[D_{s,\zeta}\hat{p}(t)|\mathcal{F}_{s}]\tfrac{\partial\tilde{\gamma
}}{\partial t}(t,s,\zeta)\nu(d\zeta)dt]ds\nonumber\\
&  =\mathbb{E}[%
{\textstyle\int_{0}^{T}}
{\textstyle\int_{t}^{T}}
{\textstyle\int_{\mathbb{R}_{0}}}
\mathbb{E}[D_{t,\zeta}\hat{p}(s)|\mathcal{F}_{t}]\tfrac{\partial\tilde{\gamma
}}{\partial s}(s,t,\zeta)\nu(d\zeta)dsdt]\nonumber\\
&  =\mathbb{E}[%
{\textstyle\int_{0}^{T}}
{\textstyle\int_{t}^{T}}
{\textstyle\int_{\mathbb{R}_{0}}}
\hat{r}(s,t,\zeta)\tfrac{\partial\tilde{\gamma}}{\partial s}(s,t,\zeta
)\nu(d\zeta)dsdt]. \label{eq2.13}%
\end{align}
Substituting $\left(  \ref{eq2.11}\right)  ,\left(  \ref{eq2.12}\right)  $ and
$\left(  \ref{eq2.13}\right)  $ combined with $\left(  \ref{eq3.3}\right)  $
in $\left(  \ref{eq2.8}\right)  $, yields%
\begin{equation}
J(u)-J(\hat{u})\leq\mathbb{E[}%
{\textstyle\int_{0}^{T}}
\{\mathcal{H}(t)-\widehat{\mathcal{H}}(t)-\tfrac{\partial\widehat{\mathcal{H}%
}}{\partial x}(t)\tilde{X}(t)\}dt].\nonumber
\end{equation}
\newline By the concavity of $\mathcal{H},$ we have%
\[
\mathcal{H}(t)-\widehat{\mathcal{H}}(t)\leq\tfrac{\partial\widehat
{\mathcal{H}}}{\partial x}(t)\tilde{X}(t)+\tfrac{\partial\widehat{\mathcal{H}%
}}{\partial u}(t)\tilde{u}(t).
\]
Hence, since $u=\hat{u}$ is $\mathbb{G}$-adapted and maximizes the conditional
Hamiltonian,
\begin{align}
&  J(u)-J(\hat{u})\leq\mathbb{E}[%
{\textstyle\int_{0}^{T}}
\tfrac{\partial\mathcal{H}}{\partial u}(t)(u(t)-\hat{u}(t))dt]\nonumber\\
&  =\mathbb{E}[%
{\textstyle\int_{0}^{T}}
\mathbb{E}[\tfrac{\partial\mathcal{H}}{\partial u}(t)|\mathcal{G}%
_{t}](u(t)-\hat{u}(t))dt]\leq0,
\end{align}
which means that $\hat{u}$ is an optimal control. \hfill$\square$ \bigskip

\subsection{A necessary maximum principle}

Suppose that a control $u\in\mathcal{A}_{\mathbb{G}}$ is optimal and that
$\beta${\normalsize $\in\mathcal{A}_{\mathbb{G}}.$ If the function }$\lambda
${\normalsize $\longmapsto J(u+\lambda\beta)$ is well-defined and
differentiable on a neighbourhood of $0$, then
\[
\tfrac{d}{d\lambda}J(u+\lambda\beta)\mid_{\lambda=0}=0.
\]
Under a set of suitable assumptions on the coefficients, we will show that%
\[
\tfrac{d}{d\lambda}J(u+\lambda\beta)\mid_{\lambda=0}=0
\]
is equivalent to%
\[
\mathbb{E}[\tfrac{\partial\mathcal{H}}{\partial u}\left(  t\right)
\mathcal{\mid G}_{t}]=0\text{ \ }P-\text{a.s. for each }t\in\lbrack0,T].
\]
}\newline The details are as follows:\newline For each given $t\in
\lbrack0,T],$ let $\eta=\eta(t)$ be a bounded $\mathcal{G}_{t}$-measurable
random variable, let $h\in\lbrack T-t,T]$ and define%

\begin{equation}
\beta(s):=\eta1_{\left[  t,t+h\right]  }(s);s\in\left[  0,T\right]  .
\label{eq4.1}%
\end{equation}
Assume that
\begin{equation}
u+\lambda\beta\in\mathcal{A}_{\mathbb{G}},
\end{equation}
for all $\beta$ and all $u\in\mathcal{A}_{\mathbb{G}}$, and all non-zero
$\lambda$ sufficiently small. Assume that the \emph{derivative process}
$Y(t)$, defined by
\begin{equation}
Y(t)=\tfrac{d}{d\lambda}X^{(u+\lambda\beta)}(t)|_{\lambda=0}, \label{1.13}%
\end{equation}
exists.\newline Then we see that%
\begin{align*}
Y(t)  &  =%
{\textstyle\int_{0}^{t}}
(\tfrac{\partial b}{\partial x}(t,s)Y(s)+\tfrac{\partial b}{\partial
u}(t,s)\beta(s))ds\\
&  +%
{\textstyle\int_{0}^{t}}
(\tfrac{\partial\sigma}{\partial x}(t,s)Y(s)+\tfrac{\partial\sigma}{\partial
u}(t,s)\beta(s))dB(s)\\
&  +%
{\textstyle\int_{0}^{t}}
{\textstyle\int_{\mathbb{R}_{0}}}
(\tfrac{\partial\gamma}{\partial x}(t,s,\zeta)Y(s)+\tfrac{\partial\gamma
}{\partial u}(t,s,\zeta)\beta(s))\tilde{N}(ds,d\zeta),
\end{align*}
and hence%

\begin{align}
dY(t)  &  =[\tfrac{\partial b}{\partial x}(t,t)Y(t)+\tfrac{\partial
b}{\partial u}(t,t)\beta(t)+%
{\textstyle\int_{0}^{t}}
(\tfrac{\partial^{2}b}{\partial t\partial x}(t,s)Y(s)+\tfrac{\partial^{2}%
b}{\partial t\partial u}(t,s)\beta(s))ds\nonumber\\
&  +%
{\textstyle\int_{0}^{t}}
(\tfrac{\partial^{2}\sigma}{\partial t\partial x}(t,s)Y(s)+\tfrac{\partial
^{2}\sigma}{\partial t\partial u}(t,s)\beta(s))dB(s)\nonumber\\
&  +%
{\textstyle\int_{0}^{t}}
{\textstyle\int_{\mathbb{R}_{0}}}
(\tfrac{\partial^{2}\gamma}{\partial t\partial x}(t,s,\zeta)Y(s)+\tfrac
{\partial^{2}\gamma}{\partial t\partial u}(t,s,\zeta)\beta(s))\tilde
{N}(ds,d\zeta)]dt\nonumber\\
&  +(\tfrac{\partial\sigma}{\partial x}(t,t)Y(t)+\tfrac{\partial\sigma
}{\partial u}(t,t)\beta(t))dB(t)\nonumber\\
&  +%
{\textstyle\int_{\mathbb{R}_{0}}}
(\tfrac{\partial\gamma}{\partial x}(t,t,\zeta)Y(t)+\tfrac{\partial\gamma
}{\partial u}(t,t,\zeta)\beta(t))\tilde{N}(dt,d\zeta). \label{1.14}%
\end{align}
We are now ready to formulate the result:

\begin{theorem}
[Necessary maximum principle]Suppose that $\hat{u}\in$ $\mathcal{A}%
_{\mathbb{G}}$ is such that, for all $\beta$ as in \eqref{eq4.1},
\begin{equation}
\tfrac{d}{d\lambda}J(\hat{u}+\lambda\beta)|_{\lambda=0}=0 \label{eq4.5}%
\end{equation}
and the corresponding solution $\hat{X}(t),(\hat{p}(t),\hat{q}(t,t),\hat
{r}(t,t,\cdot))$ of (\ref{sde}) and (\ref{p}) exists. Then,
\begin{equation}
\mathbb{E[}\tfrac{\partial\mathcal{H}}{\partial u}(t)|\mathcal{G}_{t}%
]_{u=\hat{u}(t)}=0. \label{eq4.6}%
\end{equation}
Conversely, if \eqref{eq4.6} holds, then \eqref{eq4.5} holds.
\end{theorem}

\noindent{Proof.} \quad By considering a suitable increasing family of
stopping times converging to $T$, we may assume that all the local martingales
($dB$- and $\tilde{N}$- integrals) appearing in the proof below are
martingales. We refer to the proof of Lemma 3.2 in \cite{OS} for details. For
simplicity of notation we drop the "hat" everywhere and write $u$ in stead of
$\hat{u}$, $X$ in stead of $\hat{X}$ etc in the following. Consider
\begin{equation}%
\begin{array}
[c]{l}%
\tfrac{d}{d\lambda}J(u+\lambda\beta)|_{\lambda=0}\\
=\mathbb{E[}%
{\textstyle\int_{0}^{T}}
\{\tfrac{\partial f}{\partial x}(t)Y(t)+\tfrac{\partial f}{\partial u}%
(t)\beta(t)\}dt+\tfrac{\partial g}{\partial x}(X(T))Y(T)].
\end{array}
\label{1.15}%
\end{equation}
Applying the Itô formula, we get
\begin{align*}
&
\begin{array}
[c]{l}%
\mathbb{E[}\tfrac{\partial g}{\partial x}(X(T))Y(T)]=\mathbb{E[}p(T)Y(T)]\\
=\mathbb{E[}%
{\textstyle\int_{0}^{T}}
p(t)(\tfrac{\partial b}{\partial x}(t,t)Y(t)+\tfrac{\partial b}{\partial
u}(t,t)\beta(t))dt\\
+%
{\textstyle\int_{0}^{T}}
p(t)\{%
{\textstyle\int_{0}^{t}}
(\tfrac{\partial^{2}b}{\partial t\partial x}(t,s)Y(s)+\tfrac{\partial^{2}%
b}{\partial t\partial u}(t,s)\beta(s))ds\}dt\\
+%
{\textstyle\int_{0}^{T}}
p(t)\{%
{\textstyle\int_{0}^{t}}
(\frac{\partial^{2}\sigma}{\partial t\partial x}(t,s)Y(s)+\frac{\partial
^{2}\sigma}{\partial t\partial u}(t,s)\beta(s))dB(s)\}dt\\
+%
{\textstyle\int_{0}^{T}}
p(t)\{%
{\textstyle\int_{0}^{t}}
{\textstyle\int_{\mathbb{R}_{0}}}
(\tfrac{\partial^{2}\gamma}{\partial t\partial x}(t,s,\zeta)Y(s)+\tfrac
{\partial^{2}\gamma}{\partial t\partial u}(t,s,\zeta)\beta(s))\tilde
{N}(ds,d\zeta)\}dt\\
-%
{\textstyle\int_{0}^{T}}
Y(t)\tfrac{\partial\mathcal{H}}{\partial x}(t)dt+%
{\textstyle\int_{0}^{T}}
q(t,s)(\tfrac{\partial\sigma}{\partial x}(t,t)Y(t)+\tfrac{\partial\sigma
}{\partial u}(t,t)\beta(t))dt
\end{array}
\\
&  +%
{\textstyle\int_{0}^{T}}
{\textstyle\int_{\mathbb{R}_{0}}}
r(t,s,\zeta)(\tfrac{\partial\gamma}{\partial x}(t,t,\zeta)Y(t)+\tfrac
{\partial\gamma}{\partial u}(t,t,\zeta)\beta(t))\nu(d\zeta)dt].
\end{align*}
From $\left(  \ref{eq2.12}\right)  $ and $\left(  \ref{eq2.13}\right)  $, we
have%
\begin{align*}
&
\begin{array}
[c]{l}%
\mathbb{E}\left[  p(T)Y(T)\right] \\
=\mathbb{E}[%
{\textstyle\int_{0}^{T}}
\{\tfrac{\partial b}{\partial x}(t,t)p(t)+%
{\textstyle\int_{t}^{T}}
(\tfrac{\partial^{2}b}{\partial s\partial x}(s,t)p(s)+\tfrac{\partial
^{2}\sigma}{\partial s\partial x}(s,t)q(s,t)\\
+%
{\textstyle\int_{\mathbb{R}_{0}}}
\frac{\partial^{2}\gamma}{\partial s\partial x}(s,t,\zeta)r(s,t,\zeta
)\nu(d\zeta))ds\}Y(t)dt\\
\mathbb{+}%
{\textstyle\int_{0}^{T}}
\{\tfrac{\partial b}{\partial u}(t,t)p(t)+%
{\textstyle\int_{t}^{T}}
(\frac{\partial^{2}b}{\partial s\partial u}(s,t)p(s)+\frac{\partial^{2}\sigma
}{\partial s\partial u}(s,t)q(s,t)\\
+%
{\textstyle\int_{\mathbb{R}_{0}}}
\frac{\partial^{2}\gamma}{\partial s\partial u}(s,t,\zeta)r(s,t,\zeta
)\nu(d\zeta))ds\}\beta(t)dt\\
-%
{\textstyle\int_{0}^{T}}
\frac{\partial\mathcal{H}}{\partial x}(t)Y(t)dt+%
{\textstyle\int_{0}^{T}}
(\frac{\partial\sigma}{\partial x}(t,t)Y(t)+\frac{\partial\sigma}{\partial
u}(t,t)\beta(t))q(t,t)dt
\end{array}
\\
&  +%
{\textstyle\int_{0}^{T}}
{\textstyle\int_{\mathbb{R}_{0}}}
(\tfrac{\partial\gamma}{\partial x}(t,t,\zeta)Y(t)+\tfrac{\partial\gamma
}{\partial u}(t,t,\zeta)\beta(t))r(t,t,\zeta)\nu(d\zeta)dt].
\end{align*}
Using the definition of $\mathcal{H}$ in (\ref{eq3.3}) and the definition of
$\beta$, we obtain%

\begin{equation}
\tfrac{d}{d\lambda}J(u+\lambda\beta)|_{\lambda=0}=\mathbb{E}[%
{\textstyle\int_{0}^{T}}
\tfrac{\partial\mathcal{H}}{\partial u}(s)\beta(s)ds]=\mathbb{E}{[}%
{\textstyle\int_{t}^{t+h}}
\tfrac{\partial\mathcal{H}}{\partial u}(s)ds\alpha]. \label{eq4.8}%
\end{equation}
Now suppose that
\begin{equation}
\tfrac{d}{d\lambda}J(u+\lambda\beta)|_{\lambda=0}=0. \label{eq4.9}%
\end{equation}
Differentiating the right-hand side of \eqref{eq4.8} at $h=0$, we get%
\[
\mathbb{E}[\tfrac{\partial\mathcal{H}}{\partial u}(t)\eta]=0.
\]
Since this holds for all bounded $\mathcal{G}_{t}$-measurable $\eta$, we have%
\begin{equation}
\mathbb{E[}\tfrac{\partial\mathcal{H}}{\partial u}(t)|\mathcal{G}_{t}]=0.
\label{eq4.10}%
\end{equation}
Conversely, if we assume that (\ref{eq4.10}) holds, then we obtain
(\ref{eq4.9}) by reversing the argument we used to obtain (\ref{eq4.8}).

\hfill$\square$ \bigskip

\section{Optimal consumption of a Volterra type cash flow}

Let $X^{u}(t)=X(t)$ be a given cash flow, modelled by the following stochastic
Volterra equation:
\begin{equation}%
\begin{array}
[c]{c}%
X(t)=x_{0}+%
{\textstyle\int_{0}^{t}}
[b_{0}(t,s)X(s)-u(s)]ds+%
{\textstyle\int_{0}^{t}}
\sigma_{0}(s)X(s)dB(s)\\
+%
{\textstyle\int_{0}^{t}}
{\textstyle\int_{\mathbb{R}_{0}}}
\gamma_{0}\left(  s,\zeta\right)  X(s)\tilde{N}(ds,d\zeta);\quad t\geq0,
\end{array}
\label{eq5.12}%
\end{equation}
or, in differential form,%

\begin{equation}
\left\{
\begin{array}
[c]{l}%
dX(t)=[b_{0}(t,t)X(t)-u(t)]dt+\sigma_{0}(t)X(t)dB(t)\\
+%
{\textstyle\int_{\mathbb{R}_{0}}}
\gamma_{0}\left(  t,\zeta\right)  X(t)\tilde{N}(dt,d\zeta)+[\int_{0}^{t}%
\frac{\partial b_{0}}{\partial t}(t,s)X(s)ds]dt;\quad t\geq0.\\
X(0)=x_{0}.
\end{array}
\right.  \label{eq5.13}%
\end{equation}
We see that the dynamics of $X(t)$ contains a history or memory term
represented by the $ds$-integral$.$\newline We assume that $b_{0}(t,s),$
$\sigma_{0}(s)$ and $\gamma_{0}\left(  s,\zeta\right)  $ are given 
deterministic
\color{black}
functions of $t$, $s$, and $\zeta$, with values in $\mathbb{R}$, and that
$b_{0}(t,s)$ is continuously differentiable with respect to $t$ for each $s$.
For simplicity we assume that these functions are bounded, and we assume that
there exists $\varepsilon>0$ such that $\gamma_{0}(s,\zeta)\geq-1+\varepsilon$
for all $s,\zeta$ and the initial value $x_{0}\in%
\mathbb{R}
$. We want to solve the following maximisation problem: \newline

\begin{problem}
Find $\hat{u}\in\mathcal{A_{\mathbb{G}}},$ such that
\begin{equation}
\sup_{u}J(u)=J(\hat{u}), \label{eq6.4}%
\end{equation}
where
\begin{equation}
J(u)=\mathbb{E[}\theta X(T)+%
{\textstyle\int_{0}^{T}}
\log(u(t))dt]. \label{eq5.18}%
\end{equation}

\end{problem}

Here $\theta=\theta(\omega)$ is a given $\mathcal{F}_{T}$-measurable random
variable.\newline In this case the Hamiltonian $\mathcal{H}$ gets the form%

\begin{align}
\mathcal{H}(t,x,u,\hat{p},\hat{q},\hat{r}(\cdot))  &  =\log(u)+b_{0}%
(t,t)xp-up+\sigma_{0}(t)xq\nonumber\\
&  +%
{\textstyle\int_{\mathbb{R}_{0}}}
\gamma_{0}\left(  t,\zeta\right)  xr(\zeta)\nu(d\zeta)+%
{\textstyle\int_{t}^{T}}
\tfrac{\partial b_{0}}{\partial s}(s,t)xp(s)ds. \label{eq5.16}%
\end{align}
Suppose there exists an optimal control $\hat{u}\in\mathcal{A}_{\mathbb{G}}$
for (\ref{eq5.18}) with corresponding $\hat{X},\hat{p},\hat{q},\hat{r}.$ Then,
by the optimality maximum condition we get for each $t$, that
\[
\mathbb{E[}\tfrac{\partial}{\partial u}\mathcal{H}(t,\hat{X}(t),u,\hat
{p}(t),\hat{q}(t,s),\hat{r}(t,s,\cdot))|\mathcal{G}_{t}]_{u=\hat{u}(t)}=0,
\]
i.e.,%
\[
\mathbb{E[}\frac{1}{\hat{u}(t)}\mathbb{-}\hat{p}(t)|\mathcal{G}_{t}]=0.
\]
Hence, since $\hat{u}(t)$ is $\mathbb{G}$-adapted, we get
\begin{equation}
\hat{u}(t)=\frac{1}{\mathbb{E}[\hat{p}(t)|\mathcal{G}_{t}]}. \label{^pi}%
\end{equation}
For an optimal control $\hat{u}(t)$, the corresponding adjoint equation is
reduced to the following linear BSVIE%

\begin{equation}%
\begin{array}
[c]{ll}%
\hat{p}(t) & =\theta\mathbb{+}%
{\textstyle\int_{t}^{T}}
[b_{0}(t,s)\hat{p}(s)+\sigma_{0}(s)\hat{q}(t,s)+%
{\textstyle\int_{\mathbb{R}_{0}}}
\gamma_{0}\left(  s,\zeta\right)  \hat{r}(t,s,\zeta)\nu(d\zeta)]ds\\
& -%
{\textstyle\int_{t}^{T}}
\hat{q}(t,s)dB(s)-%
{\textstyle\int_{t}^{T}}
{\textstyle\int_{\mathbb{R}_{0}}}
\hat{r}(t,s,\zeta)\tilde{N}(dt,d\zeta);\quad t\in\lbrack0,T].
\end{array}
\label{ad_e}%
\end{equation}
To find such a solution, we proceed as in
Theorem 3.1 in Hu and Øksendal \cite{ho}
\color{black}
. Define the measure $\mathbb{Q}$ by
\[
d\mathbb{Q}=M(T)dP\text{ on }\mathcal{F}_{T},
\]
where $M(t)$ satisfies the equation%
\[
\left\{
\begin{array}
[c]{ll}%
dM(t) & =M(t^{-})[\sigma_{0}(t)dB(t)+%
{\textstyle\int_{\mathbb{R}_{0}}}
\gamma_{0}(t,\zeta)\tilde{N}(dt,d\zeta)];\quad t\in\lbrack0,T],\\
M(0) & =1,
\end{array}
\right.
\]
which has a solution
\begin{align*}
M(t)  &  :=\exp(%
{\textstyle\int_{0}^{t}}
\sigma_{0}(s)dB(s)-\tfrac{1}{2}%
{\textstyle\int_{0}^{t}}
\sigma_{0}^{2}(s)ds+%
{\textstyle\int_{0}^{t}}
{\textstyle\int_{\mathbb{R}_{0}}}
\ln(1+\gamma_{0}(s,\zeta))\tilde{N}(ds,d\zeta)\\
&  +%
{\textstyle\int_{0}^{t}}
{\textstyle\int_{{\mathbb{R}}_{0}}}
\{\ln(1+\gamma_{0}(s,\zeta))-\gamma_{0}(s,\zeta)\}\nu(d\zeta)ds);\quad
t\in\lbrack0,T].
\end{align*}
Then under the measure $\mathbb{Q}$ the process
\begin{equation}
B_{\mathbb{Q}}(t):=B(t)-%
{\textstyle\int_{0}^{t}}
\sigma_{0}(s)ds\,,\quad t\in\lbrack0,T] \label{e.def_bq}%
\end{equation}
is a $\mathbb{Q}$-Brownian motion, and the random measure
\begin{equation}
\tilde{N}_{\mathbb{Q}}(dt,d\zeta):=\tilde{N}(dt,d\zeta)-\gamma_{0}(t,\zeta
)\nu(d\zeta)dt \label{e.def_nq}%
\end{equation}
is the $\mathbb{Q}$-compensated Poisson random measure of $N(\cdot,\cdot)$, in
the sense that the process%

\[
\tilde{N}_{\gamma}(t):=%
{\textstyle\int_{0}^{t}}
{\textstyle\int_{{\mathbb{R}}_{0}}}
\chi(s,\zeta)\tilde{N}_{\mathbb{Q}}(ds,d\zeta)
\]
is a local $\mathbb{Q}$-martingale, for all predictable processes
$\chi(t,\zeta)$ such that
\[%
{\textstyle\int_{0}^{T}}
{\textstyle\int_{{\mathbb{R}}_{0}}}
\gamma_{0}^{2}(t,\zeta)\chi^{2}(t,\zeta)\nu(d\zeta)dt<\infty.
\]
For all $0\leq t\leq\delta\leq T,$ define
\[
b_{0}^{(1)}(t,\delta)=b_{0}(t,\delta)\,,\quad b_{0}^{(2)}(t,\delta)=%
{\textstyle\int_{t}^{\delta}}
b_{0}(t,s)b_{0}(s,\delta)ds
\]
and inductively
\[
b_{0}^{(n)}(t,\delta)=%
{\textstyle\int_{t}^{\delta}}
b_{0}^{(n-1)}(t,s)b_{0}(s,\delta)ds\,,n=3,4,\cdots\,.
\]
Note that if $|b_{0}(t,\delta)|\leq C$ (constant) for all $t,\delta$, then by
induction on $n\in%
\mathbb{N}
$
\[
|b_{0}^{(n)}(t,\delta)|\leq\tfrac{C^{n}T^{n}}{n!},
\]
for all $t,\delta,n$. Hence,
\[
\Psi(t,\delta):=\Sigma_{n=1}^{\infty}|b_{0}^{(n)}(t,\delta)|<\infty,
\]
for all $t,\delta$. By changing of measure, we can rewrite equation
\eqref{ad_e} as
\begin{equation}
\text{ }\hat{p}(t)=\theta+%
{\textstyle\int_{t}^{T}}
b_{0}(t,s)\hat{p}(s)ds-%
{\textstyle\int_{t}^{T}}
\hat{q}(t,s)dB_{\mathbb{Q}}(s)-%
{\textstyle\int_{t}^{T}}
{\textstyle\int_{\mathbb{R}_{0}}}
\hat{r}(t,s,\zeta)\tilde{N}_{\mathbb{Q}}(dt,d\zeta);0\leq t\leq T,
\label{e.1.2}%
\end{equation}
where the processes $B_{\mathbb{Q}}$ and $\tilde{N}_{\mathbb{Q}}$ are defined
by \eqref{e.def_bq}-\eqref{e.def_nq}. Taking the conditional $\mathbb{Q}%
$-expectation on $\mathcal{F}_{t}$, we get
\begin{align}
\hat{p}(t)  &  =\mathbb{E}_{\mathbb{Q}}[\theta+%
{\textstyle\int_{t}^{T}}
b_{0}(t,s)\hat{p}(s)ds|\mathcal{F}_{t}]\nonumber\\
&  =\tilde{F}(t)+%
{\textstyle\int_{t}^{T}}
b_{0}(t,s)\mathbb{E}_{\mathbb{Q}}[\hat{p}(s)|\mathcal{F}_{t}]ds\,,\quad0\leq
t\leq T, \label{1.1.5}%
\end{align}
where
\[
\tilde{F}(s)=\mathbb{E}_{\mathbb{Q}}[\theta|\mathcal{F}_{s}\,].
\]
Fix $\delta\in\lbrack0,t]$. Taking the conditional $\mathbb{Q}$-expectation on
$\mathcal{F}_{\delta}$ of \eqref{1.1.5}, we get
\[
{\mathbb{E}}_{\mathbb{Q}}\left[  \hat{p}(t)|{\mathcal{F}}_{\delta}\right]
=\tilde{F}(\delta)+%
{\textstyle\int_{t}^{T}}
b_{0}(t,s)\mathbb{E}_{\mathbb{Q}}[\hat{p}(s)|\mathcal{F}_{\delta}%
]ds\,,\quad\delta\leq t\leq T\,.
\]
Put
\[
\tilde{p}(s)={\mathbb{E}}_{\mathbb{Q}}\left[  \hat{p}(s)|{\mathcal{F}}%
_{\delta}\right]  \,,\quad\delta\leq s\leq T\,.
\]
Then the above equation can be written as
\[
\tilde{p}(t)=\tilde{F}(\delta)+%
{\textstyle\int_{t}^{T}}
b_{0}(t,s)\tilde{p}(s)ds\,,\quad\delta\leq t\leq T\,.
\]
Substituting $\tilde{p}(s)=\tilde{F}(\delta)+\int_{s}^{T}b_{0}(s,\alpha
)\tilde{p}(\alpha)d\alpha$ in the above equation, we obtain
\begin{align}
\tilde{p}(t)  &  =\tilde{F}(\delta)+%
{\textstyle\int_{t}^{T}}
b_{0}(t,s)\{\tilde{F}(\delta)+%
{\textstyle\int_{s}^{T}}
b_{0}(s,\alpha)\tilde{p}(\alpha)d\alpha\}ds\nonumber\\
&  =\tilde{F}(\delta)+%
{\textstyle\int_{t}^{T}}
b_{0}(t,s)\tilde{F}(\delta)ds+%
{\textstyle\int_{t}^{T}}
b_{0}^{(2)}(t,\alpha)\tilde{p}(\alpha)d\alpha\,,\quad\delta\leq t\leq
T\,.\nonumber
\end{align}
Repeating this, we get by induction
\begin{align*}
\tilde{p}(t)  &  =\tilde{F}(\delta)+%
{\textstyle\sum_{n=1}^{\infty}}
{\textstyle\int_{t}^{T}}
b_{0}^{(n)}(t,\alpha)\tilde{F}(\delta)d\alpha\\
&  =\tilde{F}(\delta)+%
{\textstyle\int_{t}^{T}}
\Psi(t,\alpha)\tilde{F}(\delta)d\alpha\,.
\end{align*}
Now for $s>\delta=t$ we have ${\mathbb{E}}_{\mathbb{Q}}(\hat{p}%
(s)|{\mathcal{F}}_{t})=\tilde{p}(s)$. Hence for $s=t$ we obtain $\hat
{p}(t)=\tilde{p}(t)$, which implies that%
\begin{align}
\hat{p}(t)  &  =\tilde{F}(t)+%
{\textstyle\int_{t}^{T}}
\Psi(t,\alpha)\tilde{F}(t)d\alpha\nonumber\\
&  =\mathbb{E}_{\mathbb{Q}}[\theta+\theta%
{\textstyle\int_{t}^{T}}
\Psi(t,\alpha)d\alpha|\mathcal{F}_{t}]. \label{y_q}%
\end{align}
Moreover, by
Theorem 4.1 in Hu and Øksendal \cite{ho},
we have that $\hat{q}(t,s)$ and $\hat{r}(t,s,\zeta)$ are $C^1$ with respect to $t$ and
\color{black}
\[
\mathbb{E}_{\mathbb{Q}}[%
{\textstyle\int_{0}^{T}}
{\textstyle\int_{0}^{T}}
(\tfrac{\partial\hat{q}(t,s)}{\partial t})^{2}dsdt+%
{\textstyle\int_{0}^{T}}
{\textstyle\int_{0}^{T}}
{\textstyle\int_{\mathbb{R}_{0}}}
(\tfrac{\partial r\hat{(}t,s,\zeta)}{\partial t})^{2}\nu(d\zeta)dsdt]<\infty.
\]
Substituting the expression for $\hat{p}(t)$ in (\ref{y_q}) into the
expression of $\hat{u}(t)$ in (\ref{^pi}) and using the Bayes' rule for conditional expectation under change of
measure, we obtain the following result:

\begin{theorem}
The optimal consumption rate $\hat{u}(t)$ for Problem 5.1 is given by
\begin{align}
\hat{u}(t)  &  =\tfrac{1}{\mathbb{E}[\mathbb{E}_{\mathbb{Q}}[\theta
+\theta{\textstyle\int_{t}^{T}} \Psi(t,\alpha)d\alpha|\mathcal{F}%
_{t}]|\mathcal{G}_{t}]]} =\frac{1}{\mathbb{E}\Big[\frac{\mathbb{E}[\frac{d
\mathbb{Q}}{d\mathbb{P}}(\theta+\theta{\textstyle\int_{t}^{T}} \Psi
(t,\alpha)d\alpha) |\mathcal{F}_{t}]}{\mathbb{E}[\frac{d \mathbb{Q}%
}{d\mathbb{P}}|\mathcal{F}_{t}]}|\mathcal{G}_{t}\Big]}.
\end{align}

\end{theorem}

\begin{remark}
Here we have used that
\begin{align}
&  \mathbb{E}[\mathbb{E}_{\mathbb{Q}}[\theta+\theta{\textstyle\int_{t}^{T}}
\Psi(t,\alpha)d\alpha|\mathcal{F}_{t}]|\mathcal{G}_{t}] =\mathbb{E}%
\Big[\frac{\mathbb{E}[\frac{d \mathbb{Q}}{d\mathbb{P}}(\theta+\theta
{\textstyle\int_{t}^{T}} \Psi(t,\alpha)d\alpha) |\mathcal{F}_{t}]}%
{\mathbb{E}[\frac{d \mathbb{Q}}{d\mathbb{P}}|\mathcal{F}_{t}]}\Big |\mathcal{G}%
_{t}\Big],
\end{align}
where $\tfrac{d\mathbb{Q}}{d\mathbb{P}}=M(T)$ is the Radon-Nikodym derivative
of $\mathbb{Q}$ with respect to $\mathbb{P}$ on $\mathcal{F}_{T},$ given by
\begin{align*}
M(t)  &  :=\exp(%
{\textstyle\int_{0}^{T}}
\sigma_{0}(s)dB(s)-\tfrac{1}{2}%
{\textstyle\int_{0}^{T}}
\sigma_{0}^{2}(s)ds+%
{\textstyle\int_{0}^{T}}
{\textstyle\int_{\mathbb{R}_{0}}}
\ln(1+\gamma_{0}(s,\zeta))\tilde{N}(ds,d\zeta)\\
&  +%
{\textstyle\int_{0}^{T}}
{\textstyle\int_{{\mathbb{R}}_{0}}}
\{\ln(1+\gamma_{0}(s,\zeta))-\gamma_{0}(s,\zeta)\}\nu(d\zeta)ds);\quad
t\in\lbrack0,T].
\end{align*}

\end{remark}

\end{document}